\documentclass[12pt]{article}
\usepackage{latexsym,amssymb,upref,amsmath,amsthm, amsfonts,authblk}
\usepackage{amssymb,amsmath,amsthm, calc, graphicx}
\usepackage{epsfig}
\usepackage{breqn}
\usepackage{footnpag}
\usepackage{rotating}
\usepackage{amsfonts}
\usepackage{setspace}
\usepackage{fullpage}
\usepackage{enumitem}
\usepackage{bbold}
\usepackage{comment}
\usepackage{pgf,tikz}
\usepackage{mathrsfs}
\usetikzlibrary{arrows}
\usepackage{yhmath}
\usepackage{hyperref}
\usepackage{authblk}
\usepackage{mathrsfs}

\newtheorem{thm}{Theorem}

\newtheorem{conjecture}[thm]{Conjecture}

\newtheorem{cor}[thm]{Corollary}

\newtheorem{que}[thm]{Question}

\DeclareMathOperator{\scc}{scc}
\DeclareMathOperator{\scp}{scp}
\DeclareMathOperator{\lcc}{lcc}

\newcommand\cF{{\mathcal F}}

\newcommand{\ignore}[1]{}

\title{On Clique Coverings of Complete Multipartite Graphs}

\begin{document}
	
	\author{
		Akbar Davoodi\thanks{ 
			School of Mathematics, Institute for Research in Fundamental Sciences (IPM),
			P.O.Box: 19395-5746, Tehran, Iran, 
			E-mail: {\tt davoodi@ipm.ir}
		} \quad 
	D\'aniel Gerbner\thanks{MTA Alfr\'ed R\'enyi Institute of Mathematics, Budapest,
		Hungary, 
E-mails: {\tt gerbner@renyi.hu, vizermate@gmail.com, abhishekmethuku@gmail.com}	
}
	  \quad Abhishek Methuku$^\dagger $ \quad M\'at\'e Vizer$^\dagger $
	}
	
	\date{
		\today}

	\maketitle
	
	\begin{abstract}
		A \textit{clique covering} of a graph $G$ is a set of cliques of $G$ such that any edge of $G$ is contained in one of these cliques, and the \textit{weight} of a clique covering is the sum of the sizes of the cliques in it. The \textit{sigma clique cover number} $\scc(G)$ of a graph $G$, is defined as the smallest possible weight of a clique covering of $G$.
        
       Let $ K_t(d) $ denote the complete $ t $-partite graph with each part of size $d$. We prove that for any fixed $d \ge 2$, we have $$\lim_{t \rightarrow \infty} scc(K_t(d))= \frac{d}{2} t\log t.$$
        
        This disproves a conjecture of Davoodi, Javadi and Omoomi \cite{djo}.
		
	\end{abstract}
	
	\vspace{4mm}
	
	\section{Introduction}

	 Let $ \cal F $ be a family of not necessarily distinct sets $ S_1, S_2,\ldots, S_n $. We define the \textit{intersection graph} of $ \cal F $ in the following way: Put a vertex $ v_i $ corresponding to the set $ S_i $, for every $ 1\leq i\leq n $, and two vertices $v_i$ and $v_j$ $(i \neq j)$ are adjacent if and only if the corresponding sets intersect. In other words, the edge set of this graph is $ \{\{v_i,v_j\} : S_i\cap S_j\neq\emptyset, \ i \neq j\} $. Intersection graphs have many applications in real life problems (e.g. see \cite{mm}).
	 
	 On the other hand, given an $ n $-vertex graph $ G $ with vertex set $ \{v_1, v_2,\ldots, v_n\} $, one can find a family of sets, $ S_1, S_2,\ldots, S_n $ such that for every $ i, j $ with $ i\neq j $, $ v_i $ and $ v_j $ are adjacent in $ G $ if and only if $ S_i\cap S_j\neq\emptyset $ (for example each $S_i$ could be defined as the set of edges incident to the vertex $v_i$). This is called a \textit{set intersection representation} of $ G $. 
     In other words, a set intersection representation of $ G $ is a function $ {\cal R}:V(G) \rightarrow {\cal P}(L) $, where ${\cal P}(L)$ is the family of all subsets of a set $ L $ of labels such that for every two distinct vertices $ v_i, v_j \in V(G) $, they are adjacent if and only if $ {\cal R}(v_i)\cap {\cal R}(v_j)\neq\emptyset $.
	 .
	Note that for a fixed graph there may be several different set intersection representations and there are many related parameters that can be studied. Often the goal is to give a set intersection representation with minimum number of labels. This number is the so-called \textit{intersection number} of $ G $, denoted by $i(G)$, see for example \cite{egp}. 
    Another goal is to minimize $\max_{v\in V(G)} |{\cal R}(v)|$ over all the set intersection representations $\cal{R}$. This minimum is called the \textit{local clique cover number} of $G$ and is denoted by $\lcc(G)$ (see \cite{bdgt,jom}). 
    
In this article, we are interested in a third parameter, where the goal is to minimize  $ \sum_{v\in V(G)}|{\cal R}(v)| $  over all the set intersection representations $\cal{R}$. This minimum is called the \textit{sigma clique cover number} or the \textit{edge clique cover sum} of $G$, and is denoted by $ \scc(G) $ (see \cite{djo,k}). 
	 
Let us explain the connection between set intersection representations and clique coverings. A \textit{clique} of size $k$ is a complete graph on $k$ vertices.
A \textit{clique covering} of a graph $G$ is a set $\mathcal{C}=\{C_1,\dots, C_s\}$ of cliques of $G$ such that any edge of $G$ is contained in one of the cliques in $\mathcal{C}$. A \textit{clique partition} is a clique covering where every edge of $G$ is contained in exactly one of the cliques in $\mathcal{C}$. The \textit{clique cover number} $cc(G)$ is the smallest number of cliques in a clique covering of $G$, while the \textit{clique partition number} $cp(G)$ is the smallest number of cliques in a clique partition of $G$.
	
	Note that there is a one-to-one correspondence between the clique coverings of $ G $ and the set intersection representations of $ G $, as follows. Let $ G $ be a given graph and let $ \cal C $ be a clique covering of $G$. We consider the set of labels $ L={\cal C} $. For every vertex $ v\in V(G) $, let $\mathcal{R} (v)$ be the set of cliques containing $v$. 
    Now, if two vertices $ x $ and $ y $ are adjacent, then by definition of clique covering, there is a clique in $ \cal C $ covering the edge $xy$. This clique appears in both sets $ \mathcal{R}(x) $ and $ \mathcal{R}(y) $. This shows $i(G)\le cc(G)$. For the other direction, observe that for every label $ l\in L $ in a set intersection representation of $G$, the set of vertices of $G$ containing $l$ forms a clique in $ G $, because these vertices are pairwise adjacent by definition. This shows $i(G)\ge cc(G)$. Therefore, $i(G) = cc(G)$. (For more details see \cite{mm}.)

	In the same vein, being a clique partition is equivalent to have the condition that for every two sets $ S_i $ and $ S_j $ with $ i\neq j $, we have $ |S_i\cap S_j|\leq 1 $ in the language of set systems.
	 
	 Erd\H os, Goodman and P\'osa \cite{egp} showed that for any graph $G$ with $n$ vertices, there exists a set $ S $ with $ \lfloor\frac{n^2}{4}\rfloor $ elements and a family $ S_1, S_2, \ldots, S_n $ of (not necessarily distinct) subsets of $ S $  such that two vertices $ v_i $ and $ v_j $ are adjacent in $ G $ if and only if the sets $ S_i $ and $ S_j $ (corresponding to $ v_i $ and $ v_j $, respectively) intersect. In other words they proved that $ i(G)\leq\lfloor\frac{n^2}{4}\rfloor $.  Furthermore, $ \lfloor\frac{n^2}{4}\rfloor $ cannot be replaced by a smaller number as shown by the complete balanced bipartite graph on $n$ vertices.
	 
\vspace{3mm}
	 
	In this paper, we consider a weighted version of the above parameters. The \textit{weight} of a set of cliques $\mathcal{C}$ is the sum of the sizes of the cliques in $\mathcal{C}$. Equivalently, the weight is $ \sum_{i=1}^n|S_i| $, where $ S_1, S_2,\ldots, S_n $ is the corresponding set intersection representation of the graph. Therefore, \textit{sigma clique cover number} $\scc(G)$ is the smallest weight of a clique covering of $G$. In other words,
	\begin{align}
\scc(G)=\min_{{\cal C}: \ \text{Clique Covering}} {\sum_{c\in {\cal C}} |c|}=\min_{S_1, S_2,\ldots, S_n} \sum_{i=1}^n|S_i|,
	\end{align} 
	where the first minimum is taken over all possible clique coverings $ \cal C $ of $G$, while the latter minimum is taken over all set intersection representations $ S_1, S_2,\ldots, S_n $ of $ G $. Similarly, the \textit{sigma clique partition number}, denoted by $ \scp(G) $, is the smallest weight of a clique partition (see for example \cite{djo2}) of $G$. That is,
		\begin{align}
		\scp(G)=\min_{{\cal C}: \ \text{Clique Partition}} {\sum_{c\in {\cal C}} |c|}=\min_{S_1, S_2,\ldots, S_n} \sum_{i=1}^n|S_i|,
	\end{align} 
	where the first minimum is taken over all possible clique partitions $ \cal C $ of $G$, while  the latter minimum is taken over all set intersection representations $ S_1, S_2,\ldots, S_n $ in which  every pair of $ S_i $'s intersect in at most one element. This extra condition is equivalent to saying that every edge of $ G $ is covered by a unique clique in the clique covering. 
    
	Therefore, it is obvious that for every graph $ G $, we have $ \scc(G)\leq\scp(G) $. In the 5th Hungarian Combinatorial Colloquium, Katona and Tarj\'an raised the following conjecture.  If $ G $ is a graph on $ n $ vertices, then we have $ \scc(G)\leq\lfloor\frac{n^2}{2}\rfloor $. This conjecture was proved independently by Gy\H ori and Kostochka \cite{gk}, Kahn \cite{k} and Chung \cite{c}. Moreover in \cite{gk} and \cite{c}, it is proved that $ \scp(G)\leq\lfloor\frac{n^2}{2}\rfloor $ and the only extremal construction is the Tur\'an bipartite graph. Note that this result is a stronger statement than the above mentioned result of Erd\H{o}s, Goodman and P\'osa \cite{egp}.  Recently, this upper bound for $ \scc $ was improved for a large class of graphs.     
	Using the probabilistic method,  Davoodi, Javadi and Omoomi in \cite{djo} showed that if $ G $ is a graph on $ n $ vertices with no isolated vertices and $ \Delta(\overline{G})=d-1$, then we have
		\begin{equation}\label{eq1}
			scc(G)\leq (e^2+1)nd\left\lceil\ln\left(\frac{n-1}{d-1}\right)\right\rceil.  
		\end{equation}
In \cite{djo} they raised the following question to see if this bound is sharp.

\begin{que}
\label{question}
For positive integers $ n, d $, how large can the sigma clique cover number of an $ n-$vertex graph be, if the maximum degree of its complement is $ d-1 $?
\end{que}

	\subsection*{Complete Multipartite Graphs}
	
	In \cite{djo} Davoodi, Javadi and Omoomi investigated the sigma clique cover number of complete multipartite graphs as they were conjectured to be examples where the upper bound \eqref{eq1} for $scc$ is sharp (at least for large $n$).  They proved the following.
	
	\begin{thm}[Davoodi, Javadi and Omoomi \cite{djo}]\label{djo0}
For positive integers $n,d$ with $n\geq 2d$, let $G$ be a complete multipartite graph on $n$ vertices with at least two parts of size $d$ and the other parts of size at most $d$. Then $scc(G)\geq nd$. Moreover, if $d$ is a prime power and $n\leq d(d+1)$, then $scc(G)=scp(G)=nd$.
	\end{thm}
    
Let $ K_t(d) $ denote the complete $ t $-partite graph with each part of size $d$. We know that $ scc(K_t(d))\leq cd^2 t\log t $ for some constant $ c $, by \eqref{eq1}. On the other hand, by Theorem \ref{djo0}, we have $scc(K_t(d))=d^2t$ when $t\leq (d+1)$ and $d$ is a prime power.
	
	\begin{conjecture}[Davoodi, Javadi and Omoomi \cite{djo}]\label{firstconj}
		There exists a function $f$ and a constant $c$, such that for every positive integers $t$ and $d$, if $t\geq f(d)$, then $$scc(K_t(d))\geq cd^2t\log t.$$
	\end{conjecture}
Let us state an equivalent reformulation.
	\begin{conjecture}[Davoodi, Javadi and Omoomi \cite{djo}]\label{djo} There exists a function $ f $ and a constant $ c > 0 $ such that if
		$d\geq 2$, $ t\geq f(d) $ and $\mathcal{F}= \{(A^1_i, A^2_i,\ldots, A^d_i) \ :\  1\leq i\leq t \} $ with the property that $ A^j_i \cap A^{j'}_{i'}=\emptyset $ if and only if $ i= i' $ and $ j\neq j' $. Then we have
		\[\sum_{i,j}|A_{i}^j| \geq c d^2 t\log t. \]
	\end{conjecture}
    
 Note that the affirmative answer to this conjecture would imply that the upper bound~\eqref{eq1} for $ scc $ is best possible up to a constant factor, at least for sufficiently large $n$. However, we will disprove this conjecture.
 
 \subsection*{Qualitatively independent partitions}

    A partition of a set into $d$ classes is called a \textit{$d$-partition}. Two partitions $P$ and $P'$ of an $n$-element set
	are called \textit{qualitatively independent} if every class of $P$ has a non-empty intersection
	with every class of $P'$. This definition has an intuitive meaning in probability theory.
	Two partitions can be generated by two independent random variables if and only if the partitions are qualitatively independent.
    
    Note that if we have a family of $d$-partitions $\mathcal{F}_i= \{(A^1_i, A^2_i,\ldots, A^d_i) \}$ ($1\leq i\leq t $) such that any two are qualitatively independent, then $ A^j_i \cap A^{j'}_{i'}=\emptyset $ if and only if $ i= i' $ and $ j\neq j' $. This is the property of the families prescribed in Conjecture \ref{djo}. On the other hand, if we are given families $\mathcal{F}_i= \{(A^1_i, A^2_i,\ldots, A^d_i) \}$ ($1\leq i\leq t $)  with $ A^j_i \cap A^{j'}_{i'}=\emptyset $ if and only if $ i= i' $ and $ j\neq j' $, then let $X:=\bigcup_{i=1}^{t}\bigcup_{j=1}^d A_i^j$ and $B^{d}_i:=X\setminus \bigcup_{j=1}^{d-1} A_i^j$. Then the families $\mathcal{F}'_i= \{(A^1_i, A^2_i,\ldots, A^{d-1}_{i},B^{d}_i) \}$ ($1\leq i\leq t $) are pairwise qualitatively independent. In other words, for each $i$ we extend the partial partition given by $\cF_i$ to a full partition by adding the unused elements of $X$ to $A_i^d$.

	Gargano, K\"orner and Vaccaro \cite{gkv} studied the following.
	 Let $N(n,d)$ be the largest cardinality
	of a family of $d$-partitions of an $n$-set under the restriction that any two partitions in the family are qualitatively independent.
	
	\begin{thm}[Gargano, K\"orner and Vaccaro \cite{gkv}]\label{korner}
		For every $d$ we have $$\lim_{n\rightarrow\infty} \dfrac{1}{n}\log N(n,d)=\dfrac{2}{d}.$$
	\end{thm}
	
Note that Gargano, K\"orner and Vaccaro stated their result with $\limsup$ instead of $\lim$, but it is not hard to see that actually it holds with $\lim$.	
 
	\subsection*{Our results}
	
	We disprove Conjecture \ref{djo} by proving the following theorem:
	
	\begin{thm}\label{upper}
		
		For any $d \ge 2$  and any $\epsilon > 0$, there exist a constant $C_d < \frac{1}{2}+\epsilon$ and $t(\epsilon)$ such that for all $t \ge t(\epsilon)$ there is a family  $\mathcal{F}_t= \{(A^1_i, A^2_i,\ldots, A^d_i) \ :\  1\leq i\leq t \} $ with the property that $ A^j_i \cap A^{j'}_{i'}=\emptyset $ if and only if $ i= i' $ and $ j\neq j' $, and 
		\[ \sum_{A^j_i \in \mathcal{F}_t}{|A^j_i|} \le C_d d t\log t. \]
		
	\end{thm}
    
    \noindent 
    We also prove the following lower bound:
    
    \begin{thm}\label{lower} Let  $d,t\geq 2$ and $\mathcal{F}= \{(A^1_i, A^2_i,\ldots, A^d_i) \ :\  1\leq i\leq t \} $ such that $ A^j_i \cap A^{j'}_{i'}=\emptyset $ if and only if $ i= i' $ and $ j\neq j' $. Then we have
    \[  \frac{d}{2} t\log t \le \sum_{A^j_i \in \mathcal{F}}{|A^j_i|}. \]
    \end{thm}
	
    \noindent 
    By Theorem \ref{upper} and Theorem \ref{lower} we get:
    
	\begin{cor}\label{corol}
		
		For any $d \ge 2$ we have $$\lim_{t \rightarrow \infty} scc(K_t(d))= \frac{d}{2} t\log t.$$
		
	\end{cor}

	\section{Proofs}
    
    \subsection*{Proof of Theorem \ref{upper}}

By Theorem \ref{korner}, we have that for any $\epsilon_1 > 0$ there is an $n(\epsilon_1)$ such that for all $n \ge n(\epsilon_1)$ there is a qualitatively independent family $\mathcal{F}_n(\epsilon_1)= \{(A^1_i, A^2_i,\ldots, A^d_i) \ :\  1\leq i\leq N(n,d) \} $ on an $n$-set with the property that  $$\dfrac{2}{d} - \epsilon_1 \le \dfrac{1}{n}\log (N(n,d)).$$
Or equivalently we have ($n \ge n(\epsilon_1)$)
$$e^{ \frac{2n}{d} - \epsilon_1 n} \le N(n,d).$$
For any fixed $d$, choose $\epsilon_1 < \frac{2}{d}$ and for $n \ge n(\epsilon_1)$ let $t:=N(n,d)$. Note that $$\sum_{A^j_i \in \mathcal{F}_n(\epsilon_1)}{|A^j_i|} = tn,$$ thus we have 

	$$\sum_{A^j_i \in \mathcal{F}_n(\epsilon_1)}{|A^j_i|}\le \dfrac{1}{2-d\epsilon_1} d t  \log t.$$	
Therefore, if we choose $\epsilon_1$ small enough compared to $\epsilon$, then we are done.

\qed

    \subsection*{Proof of Theorem \ref{lower}}
    
	Let $ A=\lbrack A_{i,j}\rbrack $ be a $ t $ by $ d $ matrix, where  $A_{i,j}$ is corresponding to the set $ A_i^j $ defined in the statement of Theorem \ref{lower}. For every fixed two columns of this matrix, two sets intersect if and only if they are not in a same row of the matrix (equivalently they do not belong to a same part of the multipartite graph). Therefore, the sets corresponding to every two columns satisfies Bollob\'as's Two Families Theorem \cite{bob}. Thus, for two consecutive columns we have

\begin{align*}
	\sum_{i=1}^t \binom{|A_i^j|+|A_i^{j+1}|}{|A_{i,j|}}^{-1}
    \leq 1.
\end{align*}
To prove the lower bound we apply Bollob\'as's theorem $d$ times. First, for the first and second columns, then for the second and the third columns and finally for the last and the first columns. To simplify the notation, let $k_{i,j}:=|A_i^j|$ for $1\leq i\leq t$, $1\leq j\leq d$ and $k_{i,d+1}:=k_{i,1}$. In this way, we have the following inequality  
\begin{align*}
	\sum_{j=1}^d\sum_{i=1}^t \binom{k_{i,j}+k_{i,j+1}}{k_{i,j}}^{-1}\leq d.
\end{align*}

Let $f(m)=\binom{m}{m/2}^{-1}$ for every even integer $m$, and let $f(x)$ be its linear extension in $\mathbb{R}_{>0}$. Then $f(x)$ is convex and note that for every integers $a_i$ and $b_i$, $\binom{a_i+b_i}{a_i}^{-1}\ge \binom{a_i+b_i}{\lfloor\frac{a_i+b_i}{2} \rfloor}^{-1} \ge f(a_i+b_i)$. Thus we have

\begin{align*}
	\sum_{j=1}^d\sum_{i=1}^t 
    f(k_{i,j}+k_{i,j+1})\leq d. 
\end{align*}
On the other hand, by Jensen's inequality,

     \[f\left( \frac{2\sum_{i,j} k_{i,j}}{td} \right) \leq
     \frac{1}{td}\sum_{j=1}^d\sum_{i=1}^t f(k_{i,j}+k_{i,j+1}).\]

Hence, 
$f\left( \frac{2\sum_{i,j} k_{i,j}}{td} \right) \leq \frac{1}{t}$,
and by $\binom{n}{r}\leq 2^n$ we conclude that 
$$ \sum_{i,j} k_{i,j}\geq \frac{td}{2} \log t. $$
	
\qed
	
	\section{Concluding remarks}
    
Let us note that Question \ref{question} remains open. Recall that it was conjectured that $K_t(d)$ was the natural candidate for showing the sharpness of \eqref{eq1}. However, Theorem \ref{upper} implies that this is not the case. This suggests that perhaps \eqref{eq1} can be improved, but Theorem \ref{lower} shows that it cannot be improved by more than a factor of $d$. If on the other hand \eqref{eq1} is sharp (for $n$ large enough), then one needs to find another candidate graph $G$ with  $\Delta(\overline{G})=d-1$ and large sigma clique cover number.

	\section*{Acknowledgment}

    Research of Davoodi was supported by a grant from IPM.
	
    \vspace{2mm}
    
    \noindent 
    Research of Gerbner and Methuku was supported by the National Research, Development and Innovation Office -- NKFIH, grant K 116769.
    
    \vspace{2mm}
    
    \noindent
    Research of Vizer was supported by the National Research, Development and Innovation Office -- NKFIH, grant SNN 116095 and K116769.

\end{document}